%% file: vsp-els.tex
\documentclass{elsart}

\usepackage{latexsym, amsfonts, graphics}
\begin{document}
\newtheorem{theorem}{Theorem}[section]
\newtheorem{proposition}[theorem]{Proposition}
\newtheorem{lemma}[theorem]{Lemma}
\newtheorem{rul}[theorem]{Rule}
\newtheorem{conjecture}[theorem]{Conjecture}
\newtheorem{definition}{Definition}
\newtheorem{corollary}[theorem]{Corollary}
\newcommand{\insertfig}[1]{\begin{center}\includegraphics{#1}\end{center}}
\newenvironment{proof}{{\bf Proof:}}{\hfill\rule{2mm}{2mm}}

\begin{frontmatter}

\title{Viterbi Sequences and Polytopes}

\author{Eric H. Kuo}
\address{Department of Computer Science \\
University of California, Berkeley \\
Berkeley, CA 94720-1776}
\ead{ekuo@cs.berkeley.edu}

\begin{abstract}
A Viterbi path of length $n$ of a discrete Markov chain is a
sequence of $n+1$ states that has the greatest probability of
ocurring in the Markov chain.  We divide the space of all Markov
chains into Viterbi regions in which two Markov chains are in the same
region if they have the same set of Viterbi paths. The Viterbi paths of regions
of positive measure are called Viterbi sequences.  Our main results
are (1) each Viterbi sequence can be divided into a prefix, periodic
interior, and suffix, and (2) as $n$ increases to
infinity (and the number of states remains fixed), the number of
Viterbi regions remains bounded.  The Viterbi regions correspond to
the vertices of a Newton polytope of a polynomial whose terms are the 
probabilities of sequences of length $n$. We characterize
Viterbi sequences and polytopes for two- and three-state Markov chains.
\end{abstract}
\end{frontmatter}

\input{intro.tex}

\input{polytope.tex}

\input{bound.tex}

\input{2-state.tex}

\input{three-state.tex}

\input{conclusion.tex}

\end{document}

%% file: intro.tex
\section{Introduction}\label{sec:intro}

Many problems in computational biology have been tackled using graphical models.
A typical setup consists of $n$ observed random variables $Y_1, \ldots,
Y_n$ and $m$ hidden random variables $X_1, \ldots, X_m$. Suppose we
observe $Y_1=\sigma_1, \ldots, Y_n=\sigma_n$.  A standard
inference problem is to find the hidden assignments $h_i$ that produce
the maximum a posteriori (MAP) log probability
\[ \delta_{\sigma_1 \cdots \sigma_n} = \min_{h_1,\ldots,h_m} -\log
(\Pr[X_1=h_1,\ldots,X_m=h_m, Y_1=\sigma_1,\ldots, Y_n=\sigma_n]), \]
where the $h_i$ range over all the possible assignments for the hidden
random variables.  However, when the parameters of the graphical model
change, the hidden assignments may also vary.  The parametric
inference problem is to solve this inference problem for all model
parameters simultaneously. This problem is addressed for pairwise
sequence alignment in~\cite{GS,WEL}.

This article takes the first step in generalizing the methods to
arbitrary graphical models. We begin our investigation with Markov chains
and their Viterbi paths,
i.e. the sequence of states with the greatest probability.  We
partition the space of 
$k$-state Markov chains such that two Markov chains are in the same
region if they have the same Viterbi paths. We ask
how many such regions are
there, and how do they border each other. Before we begin to answer
these questions, let us state the problem formally.

Let $\mathcal{M}_k$ be the space of all Markov chains with $k$ states.  Then 
$\mathcal{M}_k$ is a
$(k+1)(k-1)$ dimensional space in which the parameters include the probability
distribution of the initial state and the stochastic matrix of
transition probabilities.

We number the states of a Markov chain from 0 to $k-1$.
We will consider sequences of states produced by Markov chains. The 
{\it length} of a sequence is the number of transitions in the sequence.
The {\it zeroth state} is the initial state, and the $n$th state follows the
zeroth state after $n$ transitions.

Given a Markov chain $M$, the probability of a sequence is
the product of the initial probability of the first state and all the 
transition probabilities between consecutive states.  There are $k^{n+1}$ possible
sequences of length $n$.  A {\it Viterbi path of length} $n$ is a
sequence of $n+1$ states (containing $n$ transitions) with the highest
probability.  Viterbi paths of Markov chains can be computed in polynomial
time~\cite{For,Vit}.
A Markov chain may have more than one Viterbi path of length $n$;
for instance, if 012010 is a 
Viterbi path of length 5, then 010120 must also be a Viterbi path since 
both sequences have the same initial state and the same set of transitions,
only that they appear in a different order.  Two sequences are {\it equivalent}
if their set of transitions are the same.
The Viterbi paths of a Markov chain might not be all equivalent.  Consider
the Markov chain on $k$ states that has a uniform initial distribution and
a uniform transition matrix (i.e. $p_{ij}=\frac{1}{k}$ for all states $i,j$).
Since each sequence of length $n$ has the same probability
$\frac{1}{k^{n+1}}$, every sequence is a Viterbi path for this Markov
chain. 

We define a partition $\mathcal{R}_n$ of the space $\mathcal{M}_k$
into regions such  that two Markov chains $M_1, M_2$ are in the same
region iff they have the same 
Viterbi path(s) of length $n$.  These regions will vary greatly in size.
For instance, in the partition $\mathcal{R}_2$ of $\mathcal{M}_2$, the region of Markov 
chains with the Viterbi path 00 occupies one quarter of the entire space.
Yet the Markov chain in which all the initial and transition probabilities are 
1/2 constitutes a region by itself.  In this problem, we will concern ourselves
with only the regions which have positive measure in $\mathcal{M}_k$,
and we will call  them {\it Viterbi regions}.  If $S$ is a Viterbi
path for a Viterbi region $R$, then all other possible Viterbi
paths of $R$ must be equivalent to $S$; if $S'$ were a Viterbi
path not equivalent to $S$, we would have a relation $\Pr[S]=\Pr[S']$,
making $R$ a zero-measure region of $\mathcal{M}_k$.  We call $S$ a
{\it Viterbi sequence} if it is a Viterbi path for a Viterbi
region. We will generally let one Viterbi
sequence represent an entire class of equivalent paths.

As an example, Figure~\ref{fig:k=2,n=3} shows how the subspaces $\pi_0=1$
and $\pi_1=1$ of $\mathcal{M}_2$ are partitioned into Viterbi regions
for $n=3$.  

\begin{figure}
\insertfig{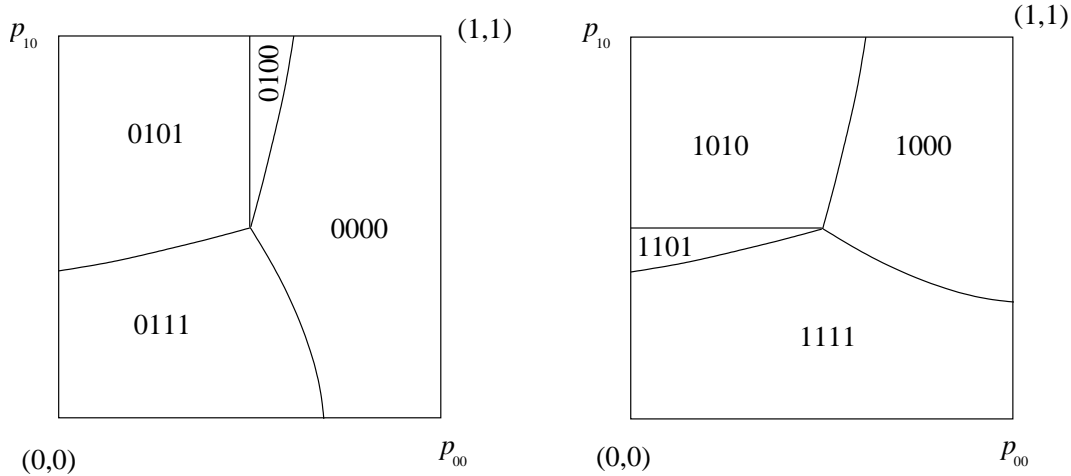}
\caption{Left: Partition of $\pi_0=1$ for $k=2$, $n=3$.  Right:
  Partition of $\pi_1=1$.}
\label{fig:k=2,n=3}
\end{figure}

In Section~\ref{sec:polytope} we use Newton polytopes to help us
enumerate all possible Viterbi regions and describe the boundary
structure within $\mathcal{M}_k$. In Section~\ref{sec:bound} we
describe the properties of Viterbi sequences and put a bound on the
number of Viterbi regions of $\mathcal{M}_k$.  
Finally, in Sections~\ref{sec:two-state} and~\ref{sec:three-state}, we
characterize Viterbi sequences for two- and three-state Markov chains.

%% file: polytope.tex
\section{Polytopes for Viterbi Sequences}\label{sec:polytope}

In analyzing the Viterbi regions of the space of all $k$-state Markov
chains, the boundary surfaces are represented as an equation between
two products of transition probabilities.  For instance, between
regions 1000 and 1111 in the space of 2-state Markov chains, the
equation is $\pi_1 p_{10}p_{00}^2 = \pi_1 p_{11}^3$.
After cancelling out $\pi_1$ on both sides, we see that the surface is
a third-degree equation.

Let us first investigate a logarithmic
model similar to a Markov chain. For each probability $p_{ij}$ or
$\pi_i$, we can define $w_{ij} = -\log p_{ij}$ as the weight of
transition $(i,j)$, and $\omega_i = -\log \pi_i$ as the weight of the
initial state $i$.  The weight of a path in the model is equal to the sum of
the weights of the transitions and initial state. A minimum weight
path is simply a path of minimum weight.

The space of logarithmic models is a subset of an even larger space
in which the weights no longer represent logarithms of
probabilities. Let $\mathcal{L}_k$ stand for the
$(k^2+k)$-dimensional space of weighting schemes.  We will often
consider only the sequences that start with state $i$; thus let
$\mathcal{L}_k^i$ be the $k^2$-dimensional subspace of $\mathcal{L}_k$,
representing the weights of the transitions (and omitting the initial
states).  The terms {\it min-weight regions} and {\it min-weight
sequences} will be defined for $\mathcal{L}_k$ and
$\mathcal{L}_k^i$ in the same way that Viterbi regions and sequences
are defined for $\mathcal{M}_k$.

Since it is not required that $\sum_{j} \exp(-w_{ij})=1$
for each state $i$, there are min-weight sequences in $\mathcal{L}_k$
that are not Viterbi sequences. These sequences will be called {\it
  pseudo-Viterbi sequences}. We will see in the next section that
001 cannot be a Viterbi sequence of a two-state Markov chain, but it
is the min-weight sequence of length 2 for weights
$(w_{00},w_{01},w_{10},w_{11}) = (3,2,4,5)$.

We can more easily visualize the min-weight regions and their boundaries
in the space $\mathcal{L}^0_k$ of linear models.  Thus instead of
considering products of probabilities, we work with sums of
logarithms.  All boundary surfaces become linear.  The boundary
between 0111 and 0000 becomes $w_{01} + 2w_{11} = 3w_{00}$.

We can formulate an alternative model for bordering min-weight
regions of $\mathcal{L}_k^0$. Instead of dividing $\mathcal{L}_k^0$ into regions
according to their min-weight sequences of length 
$n$, we shall represent each sequence as a point in a
$k^2$-dimensional space.  In this new space, each coordinate $x_{AB}$
corresponds to the number of transitions from $A$ to $B$ in the
sequence.  For example, the sequence 021021021010 is represented by
the 9-tuple $(0,1,3,4,0,0,0,3,0)$ in the case $k=3$ states.  
Note that the sequences of length $n$ occupy a $k^2-1$ dimensional 
subspace of $\mathbb{R}^{k^2}$ since the final coordinate is determined by 
the other $k^2-1$ coordinates.

Once we plot all the points corresponding to a min-weight sequence of
length $n$ into $\mathbb{R}^{k^2}$, we can examine the convex hull of
these points. This convex hull is essentially
the {\it Newton polytope} of the polynomial 
\[ \sum_{S \in \mathcal{S}_n} \Pr[S], \]
where $\mathcal{S}_n$ is the set of all possible
paths of length $n$ generated by a $k$-state Markov chain with initial
state 0, and the
probability $\Pr[S]$ is expressed in terms of transition probabilities $p_{00},
p_{01}, p_{10}$, etc.  For instance, if $\mathcal{S}_2$ is the set of
sequences of length 2 starting with state 0, then we consider the
Newton polytope of the polynomial 
$p_{00}p_{00}+p_{00}p_{01}+p_{01}p_{10}+p_{01}p_{11}$.
Figure~\ref{fig:k2poly} shows the Newton polytopes of two-state min-weight
sequences for general $n$.

Each vertex (i.e., extreme point) of the Newton polytope corresponds
to a min-weight sequence. Two vertices are joined by an edge if the
corresponding $k^2$-dimensional min-weight regions in
$\mathcal{L}_k^0$ share a boundary of dimension $k^2-1$.

The collection of log-linear cones of min-weight regions is the {\it
normal fan} of the Newton polytope.  (See~\cite[\S2.1]{St} for basic
information about Newton polytopes and their normal fans.) 
The
intersection of this normal fan with the surface $\sum_j
\exp(-w_{ij})=1$ produces a logarithmic transformation of the Viterbi
regions of $\mathcal{M}_k$.  If some min-weight region  does not
intersect $\sum_j \exp(-w_{ij})=1$, then its corresponding min-weight
sequence is also a pseudo-Viterbi sequence.

We can now describe an algorithm for enumerating min-weight sequences of
length $n$ on $k$ states.
Once we have all the min-weight sequences of length $n-1$ that start
with state 0, we derive the list of length $n-1$ min-weight sequences beginning
with state $i$ by swapping 0 and $i$ in the list of min-weight
sequences that start with 0. (In terms of coordinates, we swap
$x_{0j}$ with $x_{ij}$ and $x_{j0}$ with $x_{ji}$ for each state $j$.)
Only these sequences can be the suffix of a longer min-weight sequence.
Then to the beginning of each length $n-1$ min-weight sequence, append
state 0 (thus increasing coordinate $x_{0i}$ by 1, where $i$ was the
beginning state).  We then compute the vertices of the convex hull of
this list of coordinates, using a program such as
{\tt polymake}~\cite{GJ}.  The output is the coordinates for min-weight
sequences of length $n$ that start with state 0.

%% file: bound.tex
\section{Bound on the Number of Viterbi Sequences}\label{sec:bound}

Consider the set of Markov chains with $k$ states.  We consider the number
of possible Viterbi sequences that are of length $n$.  How does this number
grow as $n$ approaches infinity?

Somewhat surprisingly, the number of Viterbi sequences remains bounded while 
$n$ continues increase.  We show for an arbitrarily large $n$, each Viterbi
sequence can be rearranged to fit a certain blueprint. This blueprint consists
of three parts: (1) a long middle periodic section in which a sequence of at
most $k$ states is repeated, which we shall call the {\it interior}; 
(2) a short section preceding the periodic
middle, which we shall call the {\it prefix}; and (3) a section following the
middle, which we shall call the {\it suffix}.  The length of the periodic
interior should be maximized so that neither the prefix nor the suffix contains
a subsequence matching a full period of the interior; otherwise we move 
that subsequence into the interior.

\begin{lemma}\label{lem:1}
Let a subsequence of $t$ transitions begin and end with state $x$ in a
Viterbi sequence $S$.  Then if another subsequence of $t$ transitions
begin and end with a different state $y$, the set of transitions in
both subsequences must be the same.
\end{lemma}

\begin{proof}
If the two subsequences do not overlap, then $S$ must contain a subsequence
$xAxByCy$, where $A,B,C$ represent sequences of states, and $xAx$ and
$yCy$ have the same length.  Because $S$ is a Viterbi sequence, we
must have $\Pr(xAxByCy) \geq \Pr(xAxAxBy)$ and $\Pr(xAxByCy) \geq
\Pr(xByCyCy)$. Cancelling common subsequences within these
inequalities shows that $\Pr(yCy) \geq \Pr(xAx)$ and $\Pr(xAx) \geq
\Pr(yCy)$, so we conclude that $\Pr(xAx) = \Pr(yCy)$.  
But if $xAx$ and $yCy$ had different sets of transitions,
then $xByCyCy$ and $xAxAxBy$ would also be Viterbi sequences,
contradicting the definition of Viterbi sequence.

If the subsequences overlap, then $S$ contains a subsequence $xAyBxCy$. Since
$yBx$ is common to both subsequences, $xAy$ and $xCy$ have the same length.
If $xAy$ and $xCy$ have different sets of transitions, then $xAyBxAy$ and
$xCyBxCy$ are Viterbi sequences, which is also impossible.
\end{proof}

\begin{proposition}\label{prop:b2}
Let $S$ be a Viterbi sequence of length $k(mk+1)$ for some positive
integer $m$.  Then $S$ can be rearranged
to have a periodic interior that lasts for at least $m+1$ periods.
\end{proposition}

\begin{proof}
Divide $S$ into $mk+1$ sections, each $k$ transitions long.  The $k$ 
transitions are contained among $k+1$ states of $S$.  Two of these states 
must be the same, so each section contains a subsequence of transitions between
two equal states.  The length of these subsequences may be from 1 to $k$.
By pigeonhole principle, at least $m+1$ subsequences have the same length.
From Lemma~\ref{lem:1}, the set of transitions must be the same in these
$m+1$ subsequences.  We can move these subsequences so that they are adjacent
to each other.
\end{proof}

\begin{proposition}\label{prop:b3}
If a Viterbi sequence has an arbitrarily long uninterrupted periodic section 
with period $p$, then the prefix has at most $kp$ transitions, and the suffix
has at most $kp$ transitions.
\end{proposition}

\begin{proof}
Suppose there were $kp+1$ transitions in the prefix.  By the pigeonhole 
principle, at least one state from the Markov chain must appear at least 
$p+1$ times in the prefix.  Call this state $q$.  
Now classify the states in the prefix into $p$ 
classes by their distance from the start modulo $p$.  From another instance
of the pigeonhole principle, state $x$ must appear $mp$ 
transitions apart within the prefix for some integer $m$ since there are two 
states that fall into the same 
class.  The $mp$ transitions between these two $q$'s must consist of $m$ 
periods that match the interior periodic section. (Otherwise we have an
equation between the products of two sets of transition probabilities.)
However, the prefix should not contain any subsequence matching a period of
the interior.  We have a contradiction; thus the prefix can have at most $kp$
transitions.  By a similar argument, the suffix also has at most $kp$ 
transitions.
\end{proof}

However, we can prove a stronger statement about the combined lengths
of the prefix and suffix.

\begin{proposition}\label{prop:b4}
For a Viterbi sequence with periodic interior of period $p$, the
combined length of the prefix and suffix may not exceed $kp+k-2p$.
\end{proposition}

\begin{proof}
We can arrange the Viterbi sequence so that no state in the interior
appears in the prefix, but may appear in the suffix at most $p-1$
times.  For if such a state $x$ did appear $p$ times in the suffix,
then along with the last appearance of $x$ in the interior, there must
be two instances of $x$ (among those $p+1$) that are apart by a
multiple of $p$.  

Now consider a state $y$ that is not in the interior.  In the proof
for Proposition~\ref{prop:b3}, it was shown that $y$ may
appear at most $p$ times in the prefix or suffix.  If $y$ appears in
both the prefix and suffix, then we can rearrange the sequence so that
$y$ appears only once in the prefix.  (This is done by moving the
transitions between the first and last $y$ in the prefix into the
suffix.)  Since $y$ can occur at most $p$ times in the suffix, $y$
appears at most $p+1$ times in the entire sequence.

Since $p$ states are in the interior and the other $k-p$ are not in
the interior, the combined length of the prefix and suffix is bounded
by $p(p-1)+(k-p)(p+1) = kp+k-2p$.
\end{proof}

Proving these propositions allows us to conclude with the following theorem.

\begin{theorem}\label{thm:4}
The number of Viterbi sequences of length $n$ for a $k$-state Markov 
chain remains bounded as $n$ approaches infinity.
\end{theorem}

\begin{proof}
We know that each sequence must have a periodic interior.  If the period is
length $p$, then there are at most
$(p-1)!{{k}\choose{p}}=k!/((k-p)!p)$ possible periods for the
interior.  Since the prefix and suffix may have at most
$kp+k-2p$ transitions (and $kp+k-2p+1$ states) between them,
there are at most $k^{kp+k-2p+1}$ choices for the prefix and suffix.
Thus an upper bound for the
number of Viterbi sequences of length $n$ (for large $n$) is
\[ \sum_{p=1}^k k^{kp+k-2p+1} \cdot \frac{k!}{(k-p)!p}. \]
\end{proof}

All the propositions proven in this section also apply to min-weight
sequences as well. 
The next theorem demonstrates how the number of min-weight sequences does
not decrease as $n$ increases as an arithmetic sequence.

\begin{theorem}\label{thm:5}
Let $K$ be the least common multiple of the first $k$ positive
integers.  If $n>K+2k^2$, then the number of min-weight sequences of length
$n+K$ is at least the number of min-weight sequences of length $n$.
\end{theorem}

\begin{proof}
Let $S$ be a min-weight sequence of length $n$.  From
Proposition~\ref{prop:b3}, the length of the prefix and suffix of $S$ can
each be at most $k^2$ since the period $p \leq k$.  Thus the length of
the periodic interior must be at least $K$.  We will show that extending
the interior by $K$ will result in another min-weight sequence $S'$. Let
$W=\{w_{ij}\}$ be a set of weights for which $S$ is the lightest
sequence of length $n$.  Let $A$ be the set of states that do not
appear in $S$. We may assume that for $i$ or $j \in A$, $w_{ij}$ is
arbitrarily large so that $S'$ is lighter than any sequence that has a
state in $A$.  We show that for this set $W$ of weights, $S'$ is the
min-weight sequence of length $n+K$.  Suppose some other sequence
$T'$ of length $n+K>2K+k^2$ is the min-weight sequence.  The periodic
interior of $T'$ must have 
length at least $2K$.  Shorten the interior section by $K$ to form a
sequence $T$ that is the same length as $S$.  If $S'$ and $T'$ had the
same period, then subtracting the same subsequence of length $K$ from
$S'$ and $T'$ should leave $T$ lighter than $S$, which is a
contradiction with $S$ being the min-weight sequence.  Thus $S'$ and $T'$
have different periods.  Recall that every state in $T'$ (and $T$)
must also appear in $S$.  Let $P_S$ and $P_T$ be the periods of $S$
and $T$ repeated so that they are both length $K$.  If $P_S$ were
heavier than $P_T$, then we could create a sequence lighter than $S$ by
removing $P_S$ from $S$ and inserting $P_T$.  Thus $P_T$ must be the
heavier than $P_S$.  But then $S' = S+P_S$ and $T'=T+P_T$, so $S'$ is
still the lighter sequence.  Hence $S'$ is the min-weight sequence.

Since we can create a new min-weight sequence of length $n+K$ from a
sequence of length $n$, and each new sequence is distinct, the number
of min-weight sequences doesn't change.
\end{proof}

Remark: Theorem~\ref{thm:5} may be true for some values of $n$ less
than $K+2k^2$.

The previous two theorems demonstrates the eventual pattern of the number
of min-weight sequences as $n$ increases to infinity.

\begin{corollary}
Let $V_k(n)$ be the number of min-weight sequences of length $n$ on $k$
states. The sequence $\{V_k(n)\}$ eventually becomes periodic with period $K$. 
\end{corollary}

\begin{proof}
Consider the subsequence $\{V_k(n+mK)\}$.  Theorem~\ref{thm:5}
demonstrated that this subsequence is nondecreasing.  However, Theorem
\ref{thm:4} demonstrated that this sequence is bounded.  Thus each
subsequence must converge, so the sequence $\{V_k(n)\}$ must be periodic.
\end{proof}

%% file: 2-state.tex
\section{Two-State Viterbi and Min-Weight Sequences}\label{sec:two-state}

In this section we will use the properties in Section~\ref{sec:bound} to
describe the min-weight sequences for two-state linear models that
start with state 0.  The only possible periods are 00, 11, or 010.
(This notation describes the {\it transitions} in the period; thus
period 00 refers to a long run of 0's, while period 010 refers
to an alternating sequence of 0's and 1's.  Periods 010 and 101 are
equivalent.)

Period 00: Suppose state 1 exists in this min-weight
sequence. Transition 11 cannot exist, for it violates
Lemma~\ref{lem:1}. Transition 10 cannot follow since 000 and 010
cannot both be subsequences.  Thus the only min-weight sequences must
have the form $0^*$ or $0^*1$, where $i^*$ means state $i$ is repeated
as many times as desired.

Period 11: Likewise, transition 00 cannot exist. Immediately following
the initial 0 is a run of 1's.  State 0 may follow again, but then it must
be the final state: 101 and 111 cannot coexist in the same sequence.
The only possible min-weight sequences are $01^*$ and $01^*0$.

Period 010: Transition 00 or 11 may appear, but not both; if either
appears, it may appear at most once.  The only possibilities are of
the form $(01)^*, (01)^*0, (01)^*1, (01)^*10, (01)^*00$, and $(01)^*001$.
Note that exactly three of these sequences have odd length, and the
other three even.

\begin{figure}
\insertfig{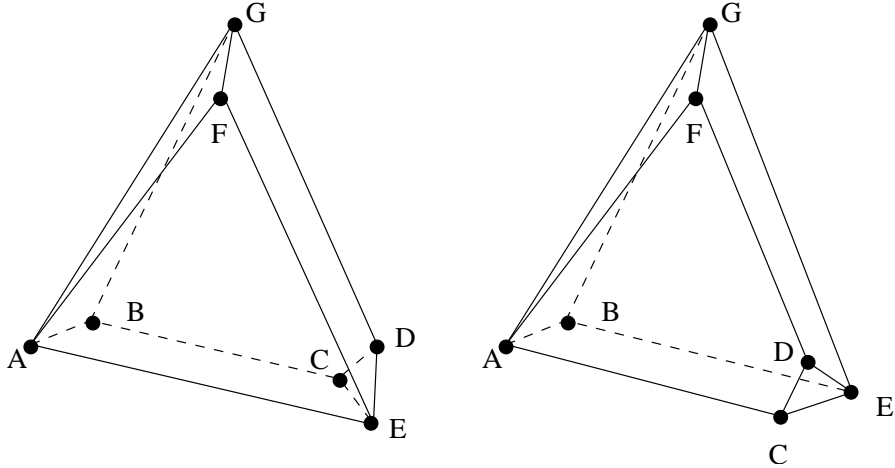}
\caption{Polytopes of two-state min-weight sequences.  The left polytope
  is for even length sequences, while the 
  right polytope is for odd length.  The vertices are labeled
  according to Table~\ref{tab:Vit2}.}
\label{fig:k2poly}
\end{figure}

Thus for $n>3$, there are at most 7 candidates for min-weight sequences.  Their
coordinates are displayed in 
Table~\ref{tab:Vit2}. If we plot these points in $\mathbb{R}^4$, their
convex hull is a three-dimensional polytope with seven vertices,
12 edges, and seven faces. Thus we can verify that all seven candidate
are indeed min-weight sequences. Different polytopes are formed for
odd and even length sequences.  These polytopes are shown in
Figure~\ref{fig:k2poly}. 

\begin{table}
\begin{center}
\begin{tabular}{lll|lll}
  $A*$ & $0^{2m+1}$ & $(2m,0,0,0)$ &  $A*$ & $0^{2m+2}$ & $(2m+1,0,0,0)$ \\
  $B$  & $0^{2m}1$  & $(2m-1,1,0,0)$ & $B$  & $0^{2m+1}1$  & $(2m,1,0,0)$ \\
  $C$ & $0(01)^m$ & $(1,m,m-1,0)$ & $C*$ & $0(01)^m0$ & $(1,m,m,0)$ \\
  $D*$  & $(01)^m1$   & $(0,m,m-1,1)$ &  $D$  & $(01)^m10$   & $(0,m,m,1)$ \\
  $E*$ & $(01)^m0$ & $(0,m,m,0)$ &   $E*$ & $(01)^{m+1}$ & $(0,m+1,m,0)$ \\
  $F$  & $01^{2m-1}0$   & $(0,1,1,2m-2)$ &   $F$  & $01^{2m}0$   & $(0,1,1,2m-1)$ \\
  $G*$ & $01^{2m}$ & $(0,1,0,2m-1)$ & $G*$ & $01^{2m+1}$ & $(0,1,0,2m)$ 
\end{tabular}
\end{center}
\caption{Left: Min-weight sequences of length $2m$.  Right:
  Min-weight sequences of length $2m+1$.  Starred
  sequences are Viterbi, unstarred are pseudo-Viterbi.}\label{tab:Vit2}
\end{table}

We could view this polytope as a facet of a greater polytope in which
we include sequences that start with state 1.  The sequences that
start with 1 form a polytope isomorphic to the one formed by sequences
that start with state 0.  These two polytopes are facets on opposite
ends of this
four-dimensional polytope, where we have introduced another coordinate
designating which state we start with.  For instance, $(01)^m0$ is now
$(0,m,m,0,0)$ and $(10)^m1$ is $(0,m,m,0,1)$ since the former starts
with state 0 and the latter starts with state 1. In this larger polytope,
an edge connects a sequence starting with state 0 to a
sequence starting with state 1 in the polytope if and only if there
exists a matrix of 
weights that produces both sequences, depending on the initial state.

Not all of the min-weight sequences are Viterbi sequences.  In fact,
for each $n>3$, three of the sequences are pseudo-Viterbi because of
the following proposition.

\begin{proposition}
No Viterbi sequence on two states can end with 001 or 110.
\end{proposition}

\begin{proof}
Suppose that 001 is a Viterbi sequence.  Then since $\Pr[001] >
\Pr[010]$, we must have $p_{00} > p_{10}$. Also, $\Pr[001]>\Pr[000]$,
so $p_{01} > p_{00}$.  Finally, $\Pr[001]>\Pr[011]$, so $p_{00} >
p_{11}$.  But then
\[ 1=p_{00}+p_{01} > p_{10}+p_{00} > p_{10}+p_{11} =1, \]
which is a contradiction.  Thus no Viterbi sequence can end with 001
(or by symmetry, 110).
\end{proof}

The remaining four min-weight sequences are actual Viterbi
sequences. Stochastic transition matrices are easily constructed to
produce these Viterbi sequences.

%% file: three-state.tex
\section{Viterbi Sequences of Three-State Markov
  Chains}\label{sec:three-state}


We have established that Viterbi sequences of 
arbitrary length must have a periodic interior as well as a prefix and a
suffix. Knowing this structure, we will begin to describe Viterbi sequences
of three-state Markov Chains.  There are eight possible periods that a
Viterbi sequence could have: 00, 11, 22, 010, 020, 121, 0120, or 0210.
We consider each period in turn, omitting periods 22, 020, and 0210
for their symmetry with periods 11, 010, and 0120, respectively.  We will also
restrict to sequences that end with state 0.

When writing the prefix of a Viterbi sequence, we will enclose the
final state in parentheses.  This final state is also the beginning of
the interior.  Similarly, the first state in the suffix is the last
state of the interior, so we will enclose it in parentheses.

To help us describe these Viterbi sequences, we will make frequent use
of the following corollary, whose proof is essentially same as in
Proposition~\ref{prop:b4}:

\begin{corollary}\label{cor:1}
If a Viterbi sequence $S$ has a periodic section of period $p$, then a
state $q$ may appear in the prefix (or suffix) at most $p$ times.  If
$q$ is in the period, then it may appear at most $p-1$ times.
\end{corollary}

Period 00: Since the final state is 0, we may assume there is no suffix.
Corollary~\ref{cor:1} tells us that since the period $p=1$, 
state 0 may not appear in the prefix, and each of the other two states
may appear at most once. So the only possible sequences
with period 00 must have the form
\[ 0^*, 10^*, 20^*, 120^*, 210^*. \]

Period 11: Only states 0 and 2 may follow a run of 1's, each at most
once.  The suffix ends with 0, so the suffix may be (1)0 or (1)20.
Similarly the prefix may contain at most one 0 and one 2.  Therefore the
Viterbi sequences with period 11 must have one of these forms:
\[ 021^*0, 201^*0, 021^*20, 201^*20 \]
or any of their suffixes (e.g. $01^*20$, $1^*0$).


Period 010: We can always rearrange a Viterbi sequence so that there
is no suffix; since we assume that the sequence ends in 0, we can
always move the subsequence 010 to the end of the sequence.

The interior may begin with either state 0 or 1.  If the interior
starts with 0, then only states 0 or 2 could
immediately precede the interior.  (State 1 would only lengthen the
interior.)  We can eliminate some impossible prefixes, using
the following facts:
\begin{enumerate}
\item The prefix may not contain subsequences 010 or 101 (for they would be
  part of the interior).
\item The prefix may not contain 020 or 121, otherwise $\Pr[020] =
  \Pr[010]$ or $\Pr[121]=\Pr[101]$.
\item States 0 and 1 may each appear at most once in the
  prefix. (Corollary~\ref{cor:1}, with $p=2$.)  
\item Wherever state 0 appears in the prefix, it must be an odd
  distance from any 0 in the interior.  (It could be even distance
  only if the only transitions in between are 01 and 10.) A similar
  statement holds for state 1.
\item State 2 may appear at most twice (Corollary~\ref{cor:1}), and if
  it does, the states must be 
  an odd distance apart.  Lemma~\ref{lem:1} says if they were an even
  distance apart, then the transitions in between would also match the
  transitions of the interior, which is impossible.  
\end{enumerate}
The only possible Viterbi sequences that satisfy all these
conditions are subsequences of
\[ 20120(10)^*, 2102(10)^*, 200(10)^*, 2100(10)^*, \]
\[ 201(10)^*, 21(10)^*, 0122(10)^*,  10220(10)^*. \]


Period 121:  We do an analysis similar to that of period 010.  For
now, we will consider Viterbi sequences that end with 210 or
2100. (The ones ending in 120 and 1200 are then derived symmetrically.)
Once again we eliminate impossible prefixes with the following
facts, many of which carry over from the case 010:
\begin{enumerate}
\item States 1 and 2 may appear only once in the prefix, and at an odd
  distance from states 1 and 2 (respectively) in the interior.
\item Subsequences 101 and 202 are impossible within the Viterbi
  sequence.
\item The subsequence 201 may not appear, otherwise we could
  create an equivalent sequence ending with 212010 (instead of 21210)
  by swapping the positions of subsequences 21 and 201.  But then we violate
  Lemma~\ref{lem:1} since 212 and 010 both appear.
\item State 0 may appear at most twice in the prefix, and at an odd
  distance if they do. 
\end{enumerate}
Finally, we need the following proposition to complete the analysis:
\begin{proposition}
A Viterbi sequence (or an equivalent sequence) cannot end in 11210, or
equivalently, 12110.
\end{proposition}

\begin{proof}
Since the sequence ends with 10, we must have $p_{10} > p_{11}$.
And since 110 is a Viterbi subsequence with higher probability than
101, we must have $p_{11}p_{10} > p_{10}p_{01}$, which means
$p_{11} > p_{01}$.  Finally, 112 has higher probability than 102,
so $p_{11}p_{12} > p_{10}p_{02}$.  Then
\[ p_{12} > \frac{p_{10}p_{02}}{p_{11}} > p_{02} \]
where we use the fact that $p_{10} > p_{11}$. Thus 
\[ 1 = p_{10}+ p_{11} + p_{12} > p_{00}+p_{01}+p_{02} = 1\]
which is a contradiction.
\end{proof}

{\it Remark:} A min-weight sequence may end with 11210 once we
remove the condition that transition probabilities starting at the
same state must sum to 1.  Thus $0(21)^*10$ and $01(21)^*10$ (as well
as $0(12)^*20$ and $02(12)^*20$) are pseudo-Viterbi sequences.  In
fact, those four sequences are the only three-state 
min-weight sequences that are pseudo-Viterbi; all other min-weight
sequences are Viterbi sequences.

Thus the possible Viterbi sequences with period 121, ending with 210,
are subsequences of
\[ 0210(21)^*0, 01(21)^*0, 00(21)^*0, 001(21)^*0, 02(21)^*0,
012(21)^*0. \]

Period 0210:  Just as in the case for periods 010 or 020, we may
assume that there is no suffix, and that the sequence ends with state
0.  Also suppose that the interior starts with state 0.  Then either
state 0 or 2 immediately precedes the interior (state 1 would extend
the interior).  
\begin{enumerate}
\item If 0 is the preceding
state, then 2 may not precede 0, for then 2002 and 2102 would be
subsequences in the same Viterbi path.
\begin{enumerate}
\item If 00 precedes the
interior, then the prefix can be extended to be 2100(0) and no further.
\item If 10 precedes the interior, only state 2 can precede 10.
(If state 1 preceded 10(0), then 11 and 00 would be in the same
sequence, violating Lemma~\ref{lem:1}. If state 0 preceded 10(0),
then 0100 and 0210 would be in the same sequence.)  The prefix cannot
be extended before 210(0) since 0210(0) would extend the interior; 1210(0) and
10210 cannot coexist; and 2210(0) violates Lemma~\ref{lem:1}.
\end{enumerate}
\item Now suppose state 2 immediately precedes the interior.  If states
12 preceded the interior, then 12(0)210 is equivalent to 121020, but 121
and 020 cannot coexist (Lemma~\ref{lem:1}).  States 22 cannot precede
the interior, for 22(0)2 and 2102 cannot coexist.  Thus only 02 can
precede the interior. The only further extensions for the prefix is
either 102(0) or 10202(0), after eliminating the following prefixes:
\begin{enumerate}
\item 002(0) and 00202(0) since 0210 is in the sequence.
\item 0102(0) since 010 and 020 cannot coexist.
\item 1102(0) and 1202(0) since 10210 is in the sequence.
\item 2102(0) and 210202(0) since 2102 can be moved to the interior.
\item 2202(0) since 2102 is in the sequence.
\end{enumerate}
\end{enumerate}

If the interior begins with states 1 or 2, the analysis is symmetric.
Thus the possible Viterbi sequences are subsequences of
\[ 2100(210)^*, 21000(210)^*, 02110(210)^*, 021110(210)^*, 102(210)^*,
1022(210)^*,\]
\[21010(210)^*, 2101010(210)^*, 120(210)^*, 12020(210)^*,
021(210)^*, 02121(210)^*. \] 

It has been established that all the possible Viterbi sequences
enumerated in this section are indeed Viterbi sequences for Viterbi
regions in $\mathcal{M}_3$.
Beginning with $n \geq 11$, the number of three-state Viterbi
sequences that start with state 0 is 89 for odd $n$ and 91 for even
$n$.  For each $n$ there are also four pseudo-Viterbi sequences:
$(01)^*002, (02)^*001, 0(12)^*20, 0(21)^*10$ for even $n$, and
$(01)^*12, (02)^*21, 0(12)^*110, 0(21)^*220$ for odd $n$.  Each Newton
polytope for three-state min-weight sequences will have 93 or 95
vertices for odd and even $n$.  Empirical evidence indicates that the
f-vectors of the Newton polytopes are periodic as $n$ increases.

%% file: conclusion.tex
\section{Conclusion}

In this article, we have discussed the parametric inference problem for Markov
chains. We have seen that the number of possible Viterbi sequences
of $k$-state Markov chains remains bounded even as the length of the
sequence increases to infinity.  We have combinatorially characterized the
structure of Viterbi sequences and enumerate them for two- and
three-state Markov chains.  One topic to be addressed in the future is
whether pseudo-Viterbi sequences exist on $k\geq 4$ states.

Viterbi sequences of Markov chains are only the tip of the iceberg of
the various statistical models whose geometries can be studied. Hidden
Markov models, parametric sequence alignment models, and binary trees
are analyzed in~\cite{PS1,PS2} with techniques such as
tropicalization and polytope propagation. A bound on explanations for
hidden Markov models is given in Corollary 8 of~\cite{PS2}.
Toric ideals of phylogenetic trees are examined in~\cite{Er}.

\section{Acknowledgments}
The author would like to thank the following persons for their
generous and helpful ideas that made this work possible: Lior Pachter,
for being a great advisor; Bernd Sturmfels and Seth Sullivant for
their discussion about Newton polytopes; and Ben Reichardt and
Nicholas Eriksson for their helpful discussions.  Data about the
polytopes was generated by {\tt polymake} ~\cite{GJ}. The author received
financial support from a National Science Foundation Graduate Research
Fellowship.